\documentclass[12pt,a4]{article}
\usepackage{setspace}
\usepackage{amssymb}
\usepackage{amsmath,amssymb}
\usepackage[T1]{fontenc}
\usepackage[bf,small,tableposition=top]{caption}
\usepackage{subfig}
\usepackage{amsthm}
\usepackage{lineno}
\usepackage[figuresright]{rotating}
\usepackage{enumerate}
\usepackage{authblk}
\usepackage[T1]{fontenc}
\usepackage[utf8]{inputenc}
\usepackage{graphics}
\usepackage{graphicx}
\usepackage{multirow}
\usepackage{lscape}
\usepackage[dvipsnames]{xcolor}
\usepackage{mathrsfs}
\usepackage{hyperref}
\usepackage{xcolor}
\usepackage{textcomp}
\usepackage{listings}
\usepackage[dvipsnames]{xcolor}
\begin{document}
	\title{A New Count Regression Model including Gauss Hypergeometric Function with an application to model demand of health services}
	\author{Deepesh Bhati$^{1\dag}$ and Ishfaq S. Ahmad$^{2}$\\
	
$^1$Department of Statistics, Central University of Rajasthan, Kishangarh, India \\
$^\dag$ Corresponding author- deepesh.bhati@curaj.ac.in \\
$^{2}$ P. G. Department of Statistics, University of Kashmir, Srinagar, India
}

\maketitle

\begin{abstract}
	In this paper, an alternative count distribution suitable for modeling over dispersed, zero vertex unimodality and monotonically decreasing data sets. 
	Though the proposed probability model includes Gauss Hypergeometric special function, it possesses simple and closed expressions for various distributional characteristics. An application to count regression modeling using a well-known data set from the National Medical Expenditure Survey is discussed by considering the length of stay in hospitalization as a dependent variable and following the proposed count model. We compare our result with the classical Negative Binomial regression model and recently proposed Uniform Poisson regression model.
\end{abstract}


\section{Introduction}
In most of the developed countries, medical expenses are primarily financed by the State or by the government. Even in developing countries, the private insurance companies through health insurance plans cover these medical expenses. Therefore modelling medical expenses is of great interest in health economics and particularly in health insurance. Length of stay (LOS) in hospitalization is one of the prominent variables which have a direct impact on medical expenses. It is also an indicator of hospital performance and a basic measure of patient’s resource consumption (Berki et al., 1984). Usually, it is a discrete random variable measured in days from the date of admission to hospital. Therefore, practitioners are always struggling in search of count distribution suitable for modeling LOS data. Typical characteristics possessed by LOS dataset are overdispersion, a high proportion of zero, right-skewed etc. Hence, the selection of appropriate count model for modeling LOS data help the decision maker to have proper budget allocation. More precisely, LOS is likely to be the outcome of covariates such as gender, patient age, nature of diseases, medical practices, and patient education status and many more. Therefore, count regression model is useful for modelling LOS in the presence of covariate. Hence the purpose of this article is two-fold, first, we propose a new count model and discussed its distributional properties and secondly, we illustrate it's application in count regression modeling.  
	
	Classical models such as Poisson, geometric distribution are being frequently used count models, however, the restriction such as equality of mean and variance in Poisson model and the successive probabilities under a geometric model decrease by a constant factor makes the model for limited use.  Hence, in consequence, many attempts have been made to develop models less restrictive than Poisson, geometric etc., and possesses all these typical characteristics.  Some useful count models studied in the statistical literature, includes the negative binomial (Bliss and Fisher, 1953), generalized Poisson (Consul, 1989) and generalized negative binomial (see Cameron and Trivedi, 1998 among others). In addition to these, various methods have been employed to develop a new class of discrete distributions like mixed Poisson method (see Karlis and Xekalaki, 2005), mixed negative Binomial (Gomez et al., 2008), discretization of a continuous family of distribution and discrete analog of continuous distribution.  Few attempts in order to derive a new class of discrete distributions involving special functions such as Hypergeometric function are available. Ong and Lee (1986) proposed a generalized non-central negative binomial distribution with probability mass function (PMF)  
\begin{equation}
p(x)= \frac{(\nu)x}{x!} p^k q^\nu (1-\theta)^{l+\nu} (1-q\theta)^{-l-\nu-x} {}_2F_1\left(-x,-l;\nu;\theta q \right), \quad x=0,1,\dots  
\end{equation} 
where $\nu,l>0, 0<p=1-q<1, 0<\theta<1$ and ${}_2F_1\left(a,b;c;z\right)=\sum_{k=0}^{\infty}\frac{(a)_k(b)_k}{(c)_k}\frac{z^k}{k!}$, is the Gauss Hypergeometric function (see Abramowitz and Stegun, 1972). Later, Gupta and Ong (2004) introduced generalized negative binomial using mixed Poisson method by mixing generalized gamma distribution for poisson parameter. However, it does not exhibit properties such as mean, variance in closed form. Later, Ghitany et al. (2002) proposed Hypergeometric negative binomial distribution with PMF 
\begin{equation*}
p(x)=\frac{\alpha^a(\alpha+1)^{p-a}(p)_x}{x!(\alpha+2)^{x+p}} {}_2F_1\left(a,x+p;p;\frac{1}{\alpha+2}\right), \quad  x=0,1,2,\dots,
\end{equation*}
suitable for over dispersed, skewed and leptokurtic datasets. Recently, Chakraborty and Ong (2015) proposed COM-Poisson type negative Binomial distribution with PMF 
\begin{equation}
p(x)= \frac{(\nu)_xp^x}{(x!)^\alpha {}_1H_{\alpha-1}(\nu;1;p)}, \quad x=0,1,\dots, 
\end{equation}
where $_1H_{\alpha-1}(\nu;1;p)=\sum_{k=0}^{\infty}\frac{(\nu)_xp^x}{(x!)^\alpha }$ is a generalized Hypergeometric function.

In all the above models, parameters to be estimated are at least three, which makes the estimation quite cumbersome.  Hence, in this article, we introduce a two-parameter count distribution considered as an alternative to negative binomial distribution and suitable for count data. Moreover, distributional properties such as expectation, variance, generating function of the proposed model are also in closed form. Other motivation of the proposed model lies in the fact that the proposed distribution is monotonically decreasing, probabilities decrease with varying rate and reduces to geometric distribution with particular parametric value. Hence the proposed distribution will be the useful contribution in applied statistics for discrete data analysis.  

For the sake of completeness, following details on Gauss Hypergeometric function are used throughout the article: \\
\noindent The Gauss Hypergeometric function is defined by
\begin{equation}
{}_{2}F_{1}(a,b;c;z)=\sum _{n=0}^{\infty }{\frac {(a)_{n}(b)_{n}}{(c)_{n}}}{\frac {z^{n}}{n!}}. 
\end{equation}
where, $|z|<1$ nd $a,b$ and $c$ is real number and $(q)_n$ is the Pochhammer symbol defined as 1 for $n=0$ and  $q(q+1)(q+2)\dots (q+n-1)$ for $n>0$. Using the identity $(a)_{n+1}=a(a+1)_n$, we have

\begin{equation}
\frac {d}{dz}\ {}_{2}F_{1}(a,b;c;z)={\frac {ab}{c}}\ {}_{2}F_{1}(a+1,b+1;c+1;z),
\end{equation}
and 
\begin{equation}
\frac {d}{db}\ {}_{2}F_{1}(a,b;c;z)=\frac{za}{c} {}_2\Theta_1^{(1)} \left( \begin{aligned} 
1,1|b,b+1,a+1 \\  b+1|2,c+1 \end{aligned} \Big|  ; z,z \right) 
\end{equation}
where \[
{}_2\Theta_1^{(1)} \left( \begin{aligned}   
a_1,a_2|b_1,b_2,b_3 \\ c_1|d_1,d_2 \end{aligned} \Big|  ; x_1,x_2 \right) =\sum_{m_1=0}^{\infty}  \sum_{m_2=0}^{\infty} \frac{(a_1)_{m_1}(a_2)_{m_2}(b_1)_{m_1}}{(c_1)_{m_1}} \frac{(b_2)_{m_1+m_2}(b_3)_{m_1+m_2}}{(d_1)_{m_1+m_2}(d_2)_{m_1+m_2}}\frac{x_1^{m_1}}{m_1!} \frac{x_2^{m_2}}{m_2!}
\]
is a Kamp\'e de F\'eriet - like function (Appell and Kamp\'e de F\'eriet, 1926). Higher order derivative of Hypergeometric function with respect to $b$ are given in Ancarani and Gasanew (2009).

The rest of the article is structured as follows: Section $2$ describes the construction of probability model and recurrence relation between successive probabilities. In section 2.1, we show the relation of the proposed model with other models. In Section $3$, methods of moment and method of maximum likelihood estimation are discussed. Maximum likelihood estimation in the presence of covariate is discussed in section 4. The application of the proposed model under regression setup is illustrated, by considering real world dataset from health sector of USA, in Section $5$. 

\section{A new count distribution}
The PMF of proposed distribution is defined by following stochastic representation
\begin{equation}  \label{sr}
\begin{aligned}
X|N\sim  &\mathcal{U}(N)\\
N|r,p  \sim &\mathcal{NB}(r,p)
\end{aligned} 
\end{equation}
\textit{where,} $\mathcal{U}(N)$ \textit{is the discrete uniform distribution over support $\left\lbrace0,1,\cdots,N\right\rbrace$ and} $\mathcal{NB}$ \textit{is negative binomial distribution with parameters $r>0$, $0<p<1$ and $x  \in \mathbb{N}_0$} and is obtained in Theorem 1. \\

\noindent \textbf{Theorem 1.} \textit{Let $X\sim \mathcal{UNB}(r,p)$ be a Uniform-negative binomial distribution as defined in (\ref{sr}) then pmf is given by}
\begin{equation} \label{pmf} 
p(x)=\frac{q^{x} p^r}{(1+x)} \binom{r+x-1}{x} {_2}F_1(1, r+x; 2+x; q), \qquad  x=0,1,\cdots,
\end{equation}
\textit{with  $r>0$, $0<p<1$, $p+q=1$ and ${_2}F_1(a,b;c;z)$ is the Gauss Hypergeometric function.}

\noindent \textit{Proof:} We prove the theorem using definition of conditional probability law 
\begin{equation*} 
\begin{aligned}
		p\left(x\right)=&\sum_{n=x}^{\infty}P\left(X|N=n\right)P\left(N=n\right) \\
		=&\sum_{n=x}^{\infty} \frac{1}{(n+1)} \binom{r+n-1}{n}p^r q^n\\
		=&\sum_{j=0}^{\infty}\frac{1}{(x+j+1)}\frac{\Gamma(r+x+j)}{(x+j)! (r-1)!} p^r q^{x+j}\\
		=& p^{r}q^{x}\sum_{j=0}^{\infty}\frac{\Gamma(r+j+x)}{ (j+x+1)!(r-1)!}  q^{j}\\
		=&\frac{1}{(1+x)} \binom{r+x-1}{x} p^r q^x{_2}F_1\left(1,r+x;2+x;q\right),
\end{aligned}
\end{equation*}
which proves the theorem. \hfill $\blacksquare$ \\

\begin{figure}[htp]
\begin{center}
	\includegraphics[scale=0.75]{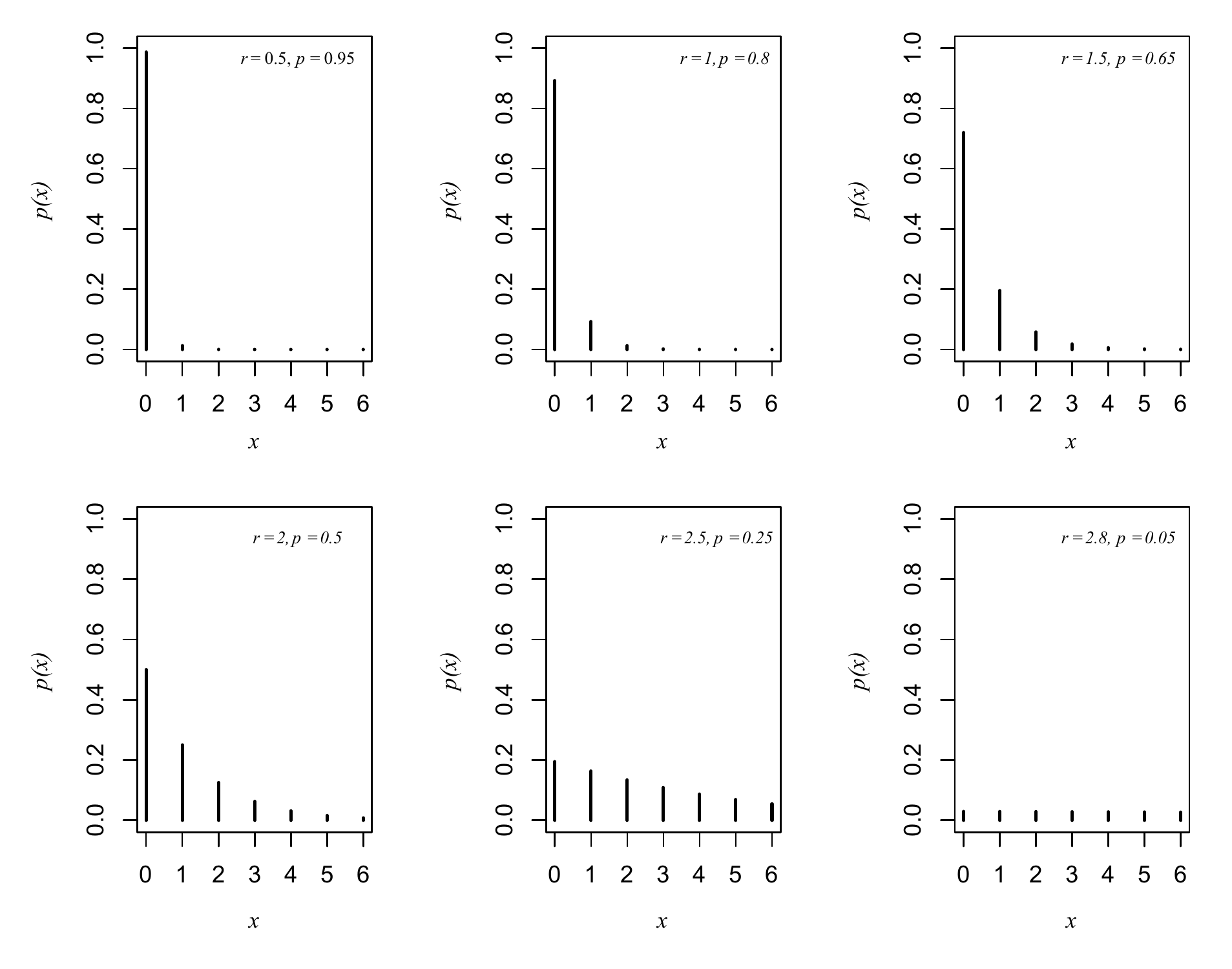}
	\captionof{figure}{PMF plot of $\mathcal{UNB}(r,p)$ distribution for different valaues of parameter $r$ and $p$.}
\end{center}
\end{figure}

Although the PMF (\ref{pmf}) contains Gauss Hypergeometric function which is a series, one can compute probabilities for different $x$ using following recurrence relation obtained by taking the ratio of probabilities at $x+1$ and $x$, and is given as
\begin{equation}
p(x+1)= q\left(\frac{r+x}{x+2}\right) \frac{_{2}F_1(1,r+x+1;3+x;q)}{_{2}F_1(1,r+x;2+x;q)} p(x).
\end{equation}
Another recurrence relation having an elegant and simple representation obtained by using the relation $\frac{\Gamma(a+j+1)}{\Gamma(a+1)}=\frac{(a)_{k+1}}{a}$ is derived as 
\begin{equation}
\begin{aligned} \nonumber 
p(x+1)&=\frac{ p^r q^{x+1}}{(2+x)} \binom{r+x}{r-1} {_2}F_1(1,r+x+1;3+x;q) \nonumber \\ \nonumber
&=\frac{p^rq^{x+1}}{(2+x)} \frac{(r+x)!}{(r-1)! (x+1)!} \sum_{j=0}^{\infty} \frac{(1)_j (r+x+1)_j}{(3+x)_j} \frac{q^j}{j!}, \\\nonumber 
&=\frac{p^{r}q^{x}}{(x+1)} \binom{r+x-1}{r-1} \sum_{t=1}^{\infty} \frac{(1)_{t} (r+x)_{t}}{(2+x)_{t}} \frac{q^{t}}{(t)!},\\ \nonumber 
&=\frac{p^{r}q^{x}}{(x+1)} \binom{r+x-1}{r-1} \left( \sum_{t=0}^{\infty} \frac{(1)_{t} (r+x)_{t}}{(2+x)_{t}} \frac{q^{t}}{(t)!}-1 \right), \\ \nonumber
&=\frac{p^{r}q^{x}}{(x+1)} \binom{r+x-1}{r-1} \left({_2}F_1(1,r+x;2+x;q)-1 \right), \\ \nonumber 
p(x+1)&=p(x)-\frac{p^{r}q^{x}}{(x+1)} \binom{r+x-1}{r-1}  \label{recc} . 
\end{aligned}
\end{equation} 

\noindent For $r>0$ and $0<p<1$, $\frac{p^{r}q^{x}}{(x+1)} \binom{r+x-1}{r-1}>0$ implies $p(x+1)<p(x)$ for all $x=0,1,\cdots$ with $p(0)=\frac{p-p^r}{(1-p)(r-1)}$. Hence it confirm that $\mathcal{UNB}(r,p)$ distribution have zero vertex unimodality. Probability plot for some parametric values are shown in Figure 1. \\

\noindent \textbf{Theorem 2.} \textit{Denoting $p_{X}{(x;r,p)}$ as PMF of $\mathcal{UNB}$ distribution with parameters  $r$ and $p$, then following recurrence relation holds true}
\begin{equation}
p_{X}{(x;r,p)}=\frac{r}{(r-1)(r+x)p}	\left(\left(2r+x-(r+x)q\right)p_{X}{(x;r+1,p)}-(r+1)p_{X}{(x;r+2,p)}\right).
\end{equation}
\noindent \textit{Proof:} Using recurrence relation of Gaussian Hypergeometric function given by 
\begin{equation*}
{_2}F_1(a,b;c;z)=\frac{2b-c+2+(a-b-1)z}{b-c+1}{_2}F_1(a,b+1;c;z)+\frac{(b+1)(z-1)}{b-c+1}{_2}F_1(a,b+2;c;z),
\end{equation*}
putting $a=1$, $b=r+x$ , $c= 2+x $ and $z=q$ in above  recurrence relation and after simplification,  we get 
\begin{equation*}
p_{X}{(x;r,p)}=\frac{r}{(r-1)(r+x)p}	\big(\left(2r+x-(r+x)q\right)p_{X}{(x;r+1,p)}-(r+1)p_{X}{(x;r+2,p)}\big),
\end{equation*}
which proves the required result. \hfill $\blacksquare$\\

\noindent \textbf{Theorem 3.} \textit{The cumulative distribution function (cdf) $F(.)$ of rv $X\sim \mathcal{UNB}(r,p)$ having PMF (\ref{pmf}) is given by}
\begin{equation*}
F_X(x)= F_Y(x)+ \frac{r+x}{x+2}\binom{r+x-1}{x} p^r q^{x+1}  {}_2F_1\left(1,r+x+1;x+3;q\right),
\end{equation*}
where $Y\sim \mathcal{NB}(r,p)$. \\
\noindent \textit{Proof:} Considering $Y\sim \mathcal{NB}(r,p)$ and using definition of distribution function, 
\begin{equation*}
\begin{aligned}
F_X(x)=P(X\le x) =&\sum_{k=0}^{x} \sum_{s=k}^{\infty}\frac{f_Y(s)}{s+1}, \\
=& \sum_{s=0}^{x} \sum_{k=0}^{s}\frac{f_Y(s)}{s+1}+\sum_{s=x+1}^{\infty} \sum_{k=0}^{x}\frac{f_Y(s)}{s+1}, \\
=& \sum_{s=0}^{x} f_Y(s) +(x+1)\sum_{s=x+1}^{\infty} \frac{f_Y(s)}{s+1}, \\
=& F_Y(x) +\sum_{s=x+1}^{\infty} \binom{r+s-1}{s} \frac{q^s p^r}{s+1},  \\
=& F_Y(x)+ \frac{x+1}{x+2}\binom{r+x}{x+1} p^r q^{x+1}  {}_2F_1\left(1,r+x+1;x+3;q\right), \\ 
F_X(x)=& F_Y(x)+ \frac{r+x}{x+2}\binom{r+x-1}{x} p^r q^{x+1}  {}_2F_1\left(1,r+x+1;x+3;q\right),
\end{aligned}
\end{equation*}
which proves the desired theorem. \hfill   $\blacksquare$\\

\subsection{Relation with existing distributions}

Using the definition of Hurwitz Lerch Zeta (HLZ) function defined by 
\begin{equation*}
\Phi \left( z,s,a\right) =\sum\limits_{k=0}^{\infty }\frac{z^{k}}{(k+a)^{s}},
\end{equation*}%
$a\in \mathbb{C}\setminus \mathbb{Z}_{0}^{-};\,s\in \mathbb{C}\,\,$when$%
\,\,|z|<1\ $and $\Re(s)>1\,\,$when$\,\,|z|=1,$ where $\mathbb{Z}^-_0$ and $\mathbb{C}$ are the set of negative integer including zero and set of complex number, respectively. In the following results, the relation of $\mathcal{UNB}(r,p)$ with $\mathcal{UG}(p)$ (Akdo$\breve{g}$an et al., 2015) is derived. \\

\noindent \textbf{Result 1:} For $r=1$, $X \stackrel{D}{=}Z-1$, where $Z\sim \mathcal{UG}(p)$ \\

\noindent \textit{Proof:} Substituting $r=1$ in (\ref{pmf}), we have

\begin{equation} 
\begin{aligned}  \nonumber 
p(x)&=\frac{ p q^{x}}{(1+x)}  {_2}F_1(1,x+1;2+x;q) \nonumber \\ \nonumber 
&=\frac{pq^{x}}{(1+x)}  \sum_{j=0}^{\infty} \frac{(1)_j (x+1)_j}{(2+x)_j} \frac{q^j}{j!}\nonumber \\ \nonumber 
&=pq^{x}  \sum_{j=0}^{\infty} \frac{q^j}{(1+x+j)}\nonumber \\ \nonumber 
&=pq^{x} \Phi(q,1,x+1) ,\quad x=0,1,\cdots,
\end{aligned} 
\end{equation}
where $q=1-p, \, \,0<p<1.$ and it is a PMF of rv $Z-1$, with $Z\sim \mathcal{UG}(p)$. \hfill $\blacksquare$\\

\noindent \textbf{Result 2:} For $r=2$, $\mathcal{UNB}(r,p)$ reduces to $\mathcal{G}(p)$. \\

\textit{Proof:} Using the fact that ${}_2F_1(1,a;a;q)=(1-q)^{-1}$, substituting $r=2$ in (\ref{pmf}), we have 
\begin{equation*}
\begin{aligned}
p(x)&=\frac{ p^{2} (1-p)^{x}}{(1+x)} \binom{x+1}{x}  {_2}F_1(1,2+x;2+x;q) \nonumber \\
&=p^2(1-p)^{x}  \sum_{j=0}^{\infty} \frac{(1)_j (2+x)_j}{(2+x)_j} \frac{q^j}{j!}\nonumber \\
&=p^{2}(1-p)^{x}  \sum_{j=0}^{\infty} q^j \\ \nonumber
&=p(1-p)^{x}  ,\quad x=0,1,\cdots,0<p<1
\end{aligned} 
\end{equation*}
which is PMF of geometric rv. Hence $\mathcal{UNB}(r,p)$ can also be viewed as generalization of geometric distribution. \hfill $\blacksquare$ \\

\noindent \textbf{Result 3:} (\textit{Limiting Distribution}) Let $r \rightarrow \infty$, $P=\frac{q}{p}\rightarrow 0$   and  $rP\rightarrow \lambda$(finite) then  $\mathcal{UNB}(r,p)$ tends to $\mathcal{UP}(\lambda)$ (see Gomez E., 2013). \\
\noindent \textit{Proof:} Assume $p=\frac{1}{Q}$ and  $q=\frac{P}{Q}$, gives  $Q-P=1$ and using condition $rP=\lambda$(finite), then (\ref{pmf}) rewritten as
\begin{equation*}
\begin{split}
p(x)=&\frac{ \left(  \frac{P}{Q}\right)  ^{x}  \left(  \frac{1}{Q} \right)  ^{r}}{(1+x)} \frac{(x+r-1) (x+r-2)...(r+1)r}{x!} \sum_{k=0}^{\infty} \frac{(1)_k (r+x)_k}{(2+x)_k} \left( \frac{P}{Q}\right) ^k\\
=&\frac{(rp)^x}{(x+1)!} \left( \frac{1}{Q}\right) ^{r+x}  \left( 1 + \frac{x-1}{r}\right) \left( \frac{1 + x-2}{r}\right) ......\left( 1 + \frac{1}{r}\right) 1 \sum_{k=0}^{\infty} \frac{(1)_k (r+x)_k }{(2+x)_k} \left( \frac{P}{1+P}\right) ^k\\
=& \frac{1}{(x+1)!} \left(  1+p\right)  ^{-r} \left( \frac{rp}{1+P}\right) ^x \sum_{k=0}^{\infty} \frac{(1)_k (r+x)_k}{(2+x)_k} \left( \frac{P}{1+P}\right) ^k\\
=&\frac{1}{(x+1)!} \left(  1+p\right)  ^{-r}  \left( \frac{\lambda}{1+ \frac{\lambda}{r}} \right) ^{x} \sum_{k=0}^{\infty} \frac{(1)_k (r+x)_k}{(2+x)_k} \left( \frac{P}{1+P}\right) ^k\\
=&\frac{1}{(1+x)!} \left(  1+ \frac{\lambda}{r}\right)  ^{-r} \left( \frac{\lambda}{1 + \frac{\lambda}{r}} \right) ^{x} \sum_{k=0}^{\infty} \frac{\Gamma(k+1)}{\Gamma(1)} \frac{\Gamma(r+x+k)}{\Gamma(r+x)}  \frac{\Gamma(2+x)}{\Gamma(2+x+k)}  \left( \frac{P}{1+P}\right) ^k\\
=& \frac{\lambda^{x}}{(x+1)!} e^{-\lambda} \sum_{k=0}^{\infty} \frac{(1)_k}{(2+x)_k} \frac{r+x+k-1}{r^{k} (r+x-1)!} \frac{\lambda^{k}}{k!},
\end{split}
\end{equation*}
\noindent as $P\rightarrow 0$,  $(1+P)^k\rightarrow 1$, hence \\
\begin{equation} \label{10}
p(x)=\frac{\lambda^{x} e^{-\lambda}}{(x+1)!} \sum_{k=0}^{\infty} \frac{(1)_k}{(2+x)_k} \lambda^{k} =\frac{\lambda^xe^{-\lambda}}{(x+1)!}{}_1F_1(1,x+2;\lambda).
\end{equation}
which is pmf of $\mathcal{UP}(\lambda)$ proposed by G\'omez (2013). \hfill $\blacksquare$

\subsection{Distribution Properties}
\subsubsection{Moment Generating Function}
The moment generating function (mgf) of $\mathcal{UNB}(r,p)$ is given by

\begin{equation*} \nonumber 
\begin{aligned} 
\mathbb{M}_X(t)=  \mathbb{E}(e^{tX})=&\sum_{x=0}^{\infty} e^{tx} \sum_{n=x}^{\infty}\frac{1}{(n+1)} \binom{r+n-1}{r-1} p^r q^n  \nonumber \\
=&\sum_{n=0}^{\infty}  \left( \sum_{x=0}^{n} \frac{e^{tx}}{n+1}\right)\binom{r+n-1}{r-1} p^r q^n \nonumber \\
=&\frac{p^r}{(e^{t}-1)} \sum_{n=0}^{\infty} \binom{r+n-1}{r-1} \left( \frac{e^{t(n+1)}-1}{n+1}\right)    q^n \nonumber \\
=&  \frac{ p^r}{e^{t}-1} \left[ \sum_{n=0}^{\infty} \binom{r+n-1}{r-1} \frac{(qe^{t})^n}{n+1} - \sum_{n=0}^{\infty} \binom{r+n-1}{r-1} \frac{q^n}{n+1}\right] \nonumber \\
=& \frac{p^r}{q(e^{t}-1)(r-1)}\left((1-p^{1-r})- e^{-t} \left(1-(1-qe^t)^{1-r}\right)\right)   \quad \text{where}  \quad t< -\log q. 
\end{aligned}
\end{equation*}

Note that using $\mathbb{M}_{X}(t)$, probability generating function of  $\mathcal{UNB}(r,p)$ distribution can also be obtained as
\begin{equation}
\begin{aligned}
\mathbb{P}_{X}(t)= &  \mathbb{E}(t^X) \nonumber\\
=& \mathbb{M}_{X}(ln(t))  \nonumber \\
=&\frac{p^r}{q(t-1)(r-1)}\left((1-p^{1-r})+t\left(1-(1-qt)^{1-r}\right)\right) 
\end{aligned}
\end{equation}

The raw moments of $\mathcal{UNB}(r,p)$ distribution can be obtained from  mgf. However, we use the stochastic representation of the $\mathcal{UNB}(r,p)$ distribution to compute the first two moments and variance. 

	\begin{equation}
		\begin{aligned}
		\mathbb{E}(X)=&\mathbb{E}_{N}(\mathbb{E}(X|N)) \nonumber\\
		=&\mathbb{E}_{N}\left(\frac{N}{2}\right) 
		=\frac{1}{2}\frac{rq}{p} \nonumber\\
		\end{aligned}
		\end{equation}
		\begin{equation}
		\begin{aligned}
		\mathbb{E}(X^2)=&\mathbb{E}_{N}(\mathbb{E
		}(X^2|N)) \nonumber\\
	=& \frac{1}{6}\left[2\mathbb{E}_{N}(N^2)+\mathbb{E}_{N}(N)\right]\nonumber\\
					=&\frac{1}{6}\left[3\frac{rq}{p}+2r(r+1)\frac{q^2}{p^2}\right]\nonumber\\
		\end{aligned}
		\end{equation}
		
\noindent Hence, the variance of $\mathcal{UNB}(r,p)$ distribution is given as
		\begin{equation}
		\begin{aligned} \nonumber 
		\mathbb{V}(X)=& \mathbb{E}(X^2)-(\mathbb{E}(X))^{2}  \nonumber \\ \nonumber 
		=&\frac{1}{6}\left(\frac{3rq}{p}+\frac{2r(r+1)q^{2}}{p^2}\right)-\frac{1}{4} \frac{r^2q^2}{p^2} \\ \nonumber
		=&\frac{rq}{12p}\left(6+\frac{4q}{p}+\frac{rq}{p}\right) 
		\end{aligned}
		\end{equation}
Further, the index of dispersion (ID) = $\frac{\mathbb{V}(X)}{	\mathbb{E}(X)} =1+\frac{4q}{6p}+\frac{rq}{6p} >1$ implies $\mathcal{UNB}(r,p)$ is always over dispersed.

\section{Parameter estimation}
			In this Section, we discuss two popular methods of estimation namely   Method of moments (MM) and Maximum Likelihood Estimation (MLE) for the estimation of the parameters for $\mathcal{UNB}(r,p)$ distribution.
			
			\subsection{Method of Moments}

			Let $X_{1},X_{2},\cdots,X_{n}$ be a sample from a population with pdf or pmf $f(x|\theta_{1},\cdots,\theta_{k})$. Moment estimators are found by equating first $k$ sample moments with corresponding $k$ population moments, and solving the resulting system of simultaneous equations. 
			Thus moments estimator  can be obtained by solving the following equation as 
			%	\textcolor{red}{text}
			\begin{equation*}
			m_{1}=\frac{r(1-p)}{2p},
			\end{equation*}
			and 
			\begin{equation*} 
			m_{2}= \frac{1}{6}\left(\frac{3r(1-p)}{p}+2r(r+1)\frac{(1-p)^2}{p^2}\right),
			\end{equation*}
			where $m_{1}$ and $m_{2}$ are first and second sample moments. Solving above system of equation for $r$ and $p$, moment estimators  obtained are 
	\begin{equation} 
\hat{r}=\frac{4m_{1}^{2}}{3(m_{2}-m_{1})-4m_{1}^{2}}.
\end{equation}
			and
			\begin{equation} 
			\hat{p}=\frac{\hat{r}}{2m_{1}+\hat{r}}
			\end{equation}

			\subsection{Maximum Likelihood Estimation}
			The method of Maximum Likelihood Estimation is, by far, the most popular technique for deriving estimators. Suppose $\underline{\mathbf{x}}=\lbrace x_{1}, x_{2},\cdots,x_{n} \rbrace$  be a random sample of size $n$ from the  $\mathcal{UNB}(r,p)$ distribution with pmf $(2)$. The likelihood function is given by 
			\begin{equation} \label{likeli} 
			L(p,r|\underline{\mathbf{x}})= \prod\limits_{i=1}^{n} p(x_i)=\prod\limits_{i=1}^{n}\frac{(1-p)^{x_i} p^r}{(1+x_i)} \binom{r+x_i-1}{x_i} {_2}F_1(1, r+x_i; 2+x_i; 1-p).
			\end{equation}
			The log-likelihood function corresponding  to $(\ref{likeli})$ is obtained as 
			\begin{equation}
			\begin{split}
			\log L(p,r|\underline{\mathbf{x}})=&\sum_{i=1}^{n}x_i \log(1-p)+\sum_{i=1}^{n} \log(p^r)-\sum_{i=1}^{n} \log(1+x_i)+\sum_{i=1}^{n} \log \Gamma(r+x_i) \\ &-n \log \Gamma(r)-\sum_{i=1}^{n} \log \Gamma(x_i) +\sum_{i=1}^{n} \log{_2}F{_1}(1,r+x_i;2+x_i;1-p).
			\end{split}
			\end{equation}
			The ML Estimates $\hat{p}$ of $p$ and $\hat{r}$ of $r$, respectively, can be obtained by solving equations 
			\begin{equation*}
			\frac{\partial \log L}{\partial p}=0, \quad \text{and} \quad  			\frac{\partial \log L}{\partial r}=0.
			\end{equation*}
where		
	\begin{equation*}
		\begin{split}
		\frac{\partial \log L}{\partial p}=& -\frac{1}{1-p}\sum_{i=1}^{n}x_i+\frac{nr}{p}+\sum_{i=1}^{n}\frac{\frac{\partial}{\partial p}
		{_2}F{_1}(1,r+x;2+x;1-p)}{{_2}F{_1}(1,r+x;2+x;1-p)}\\
	=&-\frac{1}{1-p}\sum_{i=1}^{n}x_i+\frac{nr}{p}-\sum_{i=1}^{n}\frac{(r+x_i)}{(2+x_i)}\frac{{_2}F{_1}(2,r+x_i+1;3+x_i;1-p)}{{_2}F{_1}(1,r+x_i;2+x_i;1-p)},
\end{split}
	\end{equation*}
	and 
\begin{equation*} \label{14}
\begin{aligned} 
\frac{\partial \log L}{\partial r} = n \log p&+\sum_{i=1}^{n}\psi(r+x_i)-n\psi(r) \\ &+\sum_{i=1}^{n}\frac{1}{_2F_1(1,r+x_i;2+x_i;1-p)}\frac{(1-p)}{2+x_{i}} {}_2\Theta_1^{(1)} \left( \begin{aligned} 
1,1|r+x_i,r+x_i+1,2 \\  r+x_i+1|2,x_i+2 \end{aligned} \Big|  ; 1-p,1-p \right),
\end{aligned} 
\end{equation*}
where $\psi(r)=\frac{d}{dr}\Gamma(r)$ is digamma function. As the above two equations  are not in closed form and hence cannot be solved explicitly. So we make use of iterative technique to find the ML estimates of $p$ and $r$ numerically by using \texttt{maxLik()} function in \texttt{R}.

The second order partical derivatives are given as follows 
\begin{equation*}
\begin{split}
\frac{\partial^2 \log L}{\partial p^2}=& \frac{-1}{(1-p)^2}\sum_{i=1}^{n}x_i-\frac{nr}{p^2} +2\sum_{i=1}^{n}\frac{(r+x_i)_2}{(2+x_i)_2}\frac{{_2}F{_1}(3,r+x_i+2;4+x_i;1-p)}{{_2}F{_1}(1,r+x_i;2+x_i;1-p)} \\ &-\sum_{i=1}^{n}\left(\frac{r+x_i}{2+x_i}\right)^2\left(\frac{{_2}F{_1}(2,r+x_i+1;3+x_i;1-p)}{{_2}F{_1}(1,r+x_i;2+x_i;1-p)}\right)^2
\end{split}
\end{equation*}

\begin{equation*}
\begin{split}
\frac{\partial^2 \log L}{\partial r \partial p}&=\frac{n}{p}-\sum_{i=1}^{n}\frac{1}{(2+x_i)}\frac{{_2}F{_1}(2,r+x_i+1;3+x_i;1-p)}{{_2}F{_1}(1,r+x_i;2+x_i;1-p)}\\
 -&2\sum_{i=1}^{n}\frac{r+x_i}{(2+x_i)_{2}}\frac{(1-p)}{_2F_1(1,r+x_i;2+x_i;1-p)} {}_2\Theta_1^{(1)} \left( \begin{aligned} 
1,1|r+x_i+1,r+x_i+2,3 \\  r+x_i+2|2,x_i+3 \end{aligned} \Big|  ; 1-p,1-p \right)\\
+& \sum_{i=1}^{n}\frac{r+x_i}{(2+x_i)^{2}}\frac{(1-p){}_2F_1(2,r+x_i+1;3+x_i;1-p)}{_2F_1(1,r+x_i;2+x_i;1-p)} {}_2\Theta_1^{(1)} \left( \begin{aligned} 
1,1|r+x_i,r+x_i+1,2 \\  r+x_i+1|2,x_i+2 \end{aligned} \Big|  ; 1-p,1-p \right)
\end{split}
\end{equation*}

\begin{equation*}
\begin{aligned} 
	\frac{\partial^2 \log L}{\partial r^2} = &\sum_{i=1}^{n}\psi'(r+x_i)-n\psi'(r) \\ &+\sum_{i=1}^{n}\left(\frac{1}{_2F_1(1,r+x_i;2+x_i;1-p)}\frac{(1)_{2}(1-p)^2}{(2+x_{i})_2} \cdot \right.\\ & \qquad \qquad \left. {}_2\Theta_1^{(2)} \left(\begin{aligned} 
		1,1,1|r+x_i,r+x_i+1,r+x_i+2,3 \\  r+x_i+1,r+x_i+2|3,x_i+4 \end{aligned} \Big|; 1-p,1-p,1-p \right)\right),
\end{aligned} 
\end{equation*}
where $\psi'(r)=\frac{\partial}{\partial r} \psi(r)$ is a trigamma function. \\

\noindent The expected Fisher information matrix is given as
\begin{equation*}
\mathbf{J_x}= \begin{bmatrix}
-\mathbb{E}\left(\frac{\partial^2\log L}{\partial r^2}\right)&-\mathbb{E}\left(\frac{\partial^2 \log L}{\partial r \partial p}\right) \\ 
-\mathbb{E}\left(\frac{\partial^2 \log L}{\partial p \partial r}\right)&-\mathbb{E}\left(\frac{\partial^2 \log L}{\partial p^2}\right) 
\end{bmatrix}
\end{equation*}
which can be approximated and written as
\begin{equation*} 
\mathbf{J_x} \approx \begin{bmatrix}
J_{rr} &J_{rp} \\ 
J_{pr}&J_{pp} 
\end{bmatrix}=\begin{bmatrix}
\frac{\partial^2 \log L}{\partial r^2}\big|_{\hat{r},\hat{p}}&\frac{\partial^2 \log L}{\partial r \partial p}\big|_{\hat{r},\hat{p}} \\
\\ 
\frac{\partial^2\log L}{\partial p \partial r}\big|_{\hat{r},\hat{p}}&\frac{\partial^2\log L}{\partial p^2}\big|_{\hat{r},\hat{p}} 
\end{bmatrix}
\end{equation*}
where $\hat{r}$ and $\hat{p}$ are the maximum likelihood estimators of $r$ and $p$ respectively. Hence, when $n$ is large and under some mild regularity conditions, we have that
\[
\sqrt{n}
\left( {\begin{array}{cc}
	r-\hat{r} \\
	p-\hat{p} \\
	\end{array} } \right)
\stackrel{a} \sim  N_{2}\left( \left({\begin{array}{cc}
	0 \\
	0 \end{array}}\right),\mathbf{J_x}^{-1}
\right),\]
where $"\stackrel{a}{\sim}"$ means approximately distributed, and $\mathbf{J_x}^{-1} $ is the inverse of $\mathbf{J_x}$. The above asymptotic normal distribution can be used to construct approximate confidence intervals for the parameters; that is, we have the asymptotic confidence intervals $\hat{r}\mp z_{1-\alpha/2}se(\hat{r})$ and $\hat{p}\mp z_{1-\alpha/2}se(\hat{p})$ for $r$ and $p$ respectively. Here, $se$ is the square root of the diagonal element of $\mathbf{J_x}^{-1}$ corresponding to each parameter (i.e., the asymptotic standard error), and $z_{(1-\alpha/2)}$ denotes the $(1-\alpha/2)$ quantile of standard normal distribution.

\subsubsection{Likelihood ratio test}
As $\mathcal{UNB}(r,p)$ reduces to $\mathcal{G}(p)$ when $r=2$ (see Result 2). These two distributions are nested therefore we can employ the likelihood ratio test criterion to test the following hypothesis:
\begin{equation*}
\begin{split}
\mathscr{H}_0: & \, \, r=2, \quad \text{sample is from geometric with parameter}  \,\,  p \\
\mathscr{H}_1: & \, \, r\ne 2 \quad \text {sample is from $\mathscr{UNB}$ with parameters $r,p$}
\end{split} 
\end{equation*}   
Writing $\Xi=(r,p)$  the parametric space, the likelihood ratio test statistic is given by
\[ LR= -2 \frac{l(\widehat{\Xi}^*,x)}{l(\widehat{\Xi},x)},\]
where $\widehat{\Xi}^*$ is the restricted ML estimates under the null hypothesis $\mathscr{H}_0$ and $\widehat{\Xi}$ is the unrestricted ML estimate under alternative hypothesis $\mathscr{H}_1$. Under the null hypothesis $\mathscr{H}_0$ , LR follows a chi-square distribution with one degree of freedom (df). Hence, at 5\% level of significance the two sided critical region for this test is given by $\lbrace LR : LR < 0.00098  \cap LR > 5.02\rbrace$. Thus there will be no evidence against the null hypothesis if $0.00098 < LR <5.02$, otherwise $\mathscr{H}_0$ will be rejected.
			
\section{Regression including Covariates}
In count regression modeling, Poisson regression has been used frequently. However, as indicated in beginning, if the response variation is greater than the response mean then Poisson regression model is not suitable. Such situation is very common in data such as insurance claim, health data etc. In these cases, Generalized Poisson  regression model (see Consul and Femoye, 1992), Negative binomial regression model (see Hilbe, 2007), Generalized negative Binomial, Poisson-Inverse Gaussian, Uniform Poisson (G\'omez, 2010) etc. In this section, we will discuss regression modeling by considering Uniform negative binomial distribution for response variable.  
			
Let $X$ be the response variable and $\mathbf{y}$ be associated $s\times 1$ vector of covariates. We consider that the response variable $X$ follow $\mathcal{UNB}$ distribution with mean $\mu(\mathbf{x})$. Furthermore, the mean of response variable linked with the explanatory variables by log linear form i.e. $\mu_{i}=\exp(\beta \mathbf{y_{i}}^{\top})$ where $\mathbf{\beta}=(\beta_{0},\beta_{1},\cdots,\beta_{s})$ and $\mathbf{y_i}=(1,y_{1i},\cdots,y_{si})$. By replacing $p$ with $\frac{r}{2\mu_{i}+r}$ we obtain the re-parametrize PMF as 
\begin{equation}
 f_X(x_{i})=\frac{\Gamma(r+x_{i})}{\Gamma(r)\Gamma(x_{i}+2)}\left(\frac{2\mu_{i}}{2\mu_{i}+r}\right)^{x_{i}}\left(\frac{r}{2\mu_{i}+r}\right)^{r}{_2}F_1\left(1, r+x_{i}; 2+x_{i};\frac{2\mu_{i}}{2\mu_{i}+r} \right),
\end{equation}

and the corresponding log-likelihood equation is given as  
	\begin{eqnarray} \label{llc}
	l &=& \sum_{i=1}^{n}\log \Gamma(r+x_{i})-n\log \Gamma(r)-\sum_{i=1}^{n}\log \Gamma(x_{i}+1) +\log 2 \sum_{i=1}^{n}x_{i}
	\nonumber\\
	&&+\sum_{i=1}^{n}x_{i}\left(e^{{\mathbf{\beta}} \mathbf{y_{i}}^{\top}}\right)+nr\log r-\sum_{i=1}^{n}\left(x_{i} +r\right)\log \left(2e^{{\mathbf{\beta}} \mathbf{y_{i}}^{\top}}+r\right)
	\nonumber\\
	&&-\sum_{i=1}^{n}\log (x_{i}+1)+\sum_{i=1}^{n}\log {_2}F_1\left(1, r+x_{i}; 2+x_{i};\frac{2e^{{\mathbf{\beta}} \mathbf{y_{i}}^{\top}}}{2e^{{\mathbf{\beta}} \mathbf{y_{i}}^{\top}}+r} \right) 
		\end{eqnarray}
		The normal equations with respect to parameters $(r,\underline{\mathbf{\beta}})$ are 
		\begin{eqnarray}
		\frac{\partial}{\partial r}l&=& n\log r+n+\sum_{i=1}^{n}\frac{(x_{i}+r)}{\left(2e^{\underline{\mathbf{\beta}}y_{i}^{T}}+r\right)}-\sum_{i=1}^{n}\log\left(2e^{\underline{\mathbf{\beta}}y_{i}^{T}}+r\right) 
		\nonumber\\
		&&+\sum_{i=1}^{n}\frac{1}{\Gamma(r+x_{i})}\left(\frac{\partial}{\partial r}\Gamma(r+x_{i})\right)-\frac{n}{\Gamma(r)}\left(\frac{\partial}{\partial r}\Gamma(r)\right) \nonumber \\
		&&+\sum_{i=1}^{n}\frac{\partial}{\partial r}\log {_2}F_1\left(1, r+x_{i}; 2+x_{i};\frac{2e^{{\mathbf{\beta}} \mathbf{y_{i}}^{\top}}}{2e^{{\mathbf{\beta}} \mathbf{y_{i}}^{\top}}+r} \right), \nonumber
		\end{eqnarray}
			\begin{eqnarray}
		\frac{\partial}{\partial \beta_{k}}l &=&\sum_{i=1}^{n}x_{i}y_{ki}-\sum_{i=1}^{n}\left(\frac{x_{i}+r}{2e^{\underline{\mathbf{\beta}}y_{i}^{\top}}+r}\right)2y_{ki}
		\nonumber \\
		&&	+\sum_{i=1}^{n}\frac{\partial}{\partial \beta_{k}}\log {_2}F_1\left(1, r+x_{i}; 2+x_{i};\frac{2e^{\underline{\beta}y_{i}^{\top}}}{2e^{\underline{\beta}y_{i}^{\top}}+r} \right), \quad k=0,1,\cdots,s\nonumber
		\end{eqnarray}
		\par 
Above $s+2$ normal equations are not in closed form  and can not be solved explicitly. Therefore by employing \texttt{FindMaxima()} function available in \texttt{Mathematica}, we obtained the numerical solutions. The Mathematica code can be made available on request to interested readers.

\section{Numerical Illustration}
We illustrate the application of $\mathcal{UNB}$ regression by considering data obtained from the National Medical Expenditure Survey (NMES) which was
conducted in 1987 and 1988. A subsample of size 4406 individuals above age 65 covered under medicare scheme is considered. This scheme provides substantial protection against health care costs. The dataset is available in the Journal of Applied Econometrics 1997 Data Archive at	\url{http://www.econ.queensu.ca/jae/1997-v12.3/deb-trivedi/}.      
The data originally used by Deb and Trivedi (1997) for analyzing various measures of health care utilizations. We use the variable \textit{length of stay after hospital admission} (HOSP) as the response variable and analyze it in the presence of 10 explanatory variables summarized in Table 1. The high degree of over dispersion (1.875) and proportion of zero (80\%) of HOSP variable leading to acceptance of the proposed model.

\begin{table}[htp] 
	\centering 
	\caption{Definition of explanatory variables and their summary statistics}
	\begin{tabular}{llcc} \hline  \hline  
		Explanatory variable &	Description	& mean	& stdev \\ \hline 
		EXCELHLTH	& Self-perceived health status,	& 0.078	& 0.268 \\
		& excellent =1, else = 0 & & \\
		POORHLTH	& Self-perceived health status, & 0.126	& 0.332 \\
		& poor=1, else=0	& & \\
		
		NUMCHRON	& Number of chronic conditions &	1.542	& 1.35 \\
		AGE	& Age of patient	& 7.402 &	0.633 \\
		MALE &	Gender; male = 1, else =0 &	0.404 &	0.491 \\
		MARRIED	& Marital status of patient, & 0.546	& 0.498 \\
		
		&Married =1 else =0	&&\\
		FAMINC	& Equals family income in \$10,000	& 2.527	& 2.925 \\
		EMPLOYED	& employment status, employed =1 els=0	& 0.103 &	0.304 \\
		PRINVINS &	equals 1 if the person is covered by  &	0.776	& 0.417 \\
		& private health insurance& & \\
		MEDICAID	& equals 1 if the person is covered by Medicaid	&0.091	&0.288 \\ \hline \hline 
	\end{tabular}
\end{table}

As the dependent variable is over dispersed, previously proposed models such as Uniform-Poisson, Negative binomial are suitable choice for regression modeling. The log-likelihood value obtained for HOSP data without covariate for Uniform-Poisson, Negative binomial and the proposed Uniform-negative binomial are -3193.28, -3009.62 and -3008.32 respectively. The highest value of log likelihood (LL) for proposed $\mathcal{UNB}$ distribution gives evidence of its superiority over other models. Overall distribution of HOSP is presented in Figure 2. Further the gender-wise and health status-wise of HOSP variable presented in Table 2. It is clear from both the tables that proportion of female admitted are more than male. Though the proportion of patient with poor health condition are 12\%, but there mean stay and variance indicate larger index of dispersion (2.10). The positive correlation coefficient (0.233) between HOSP and number of chronic conditions indicates the dependence.

\begin{figure}[htp]
	\begin{center}
		\includegraphics[scale=0.75]{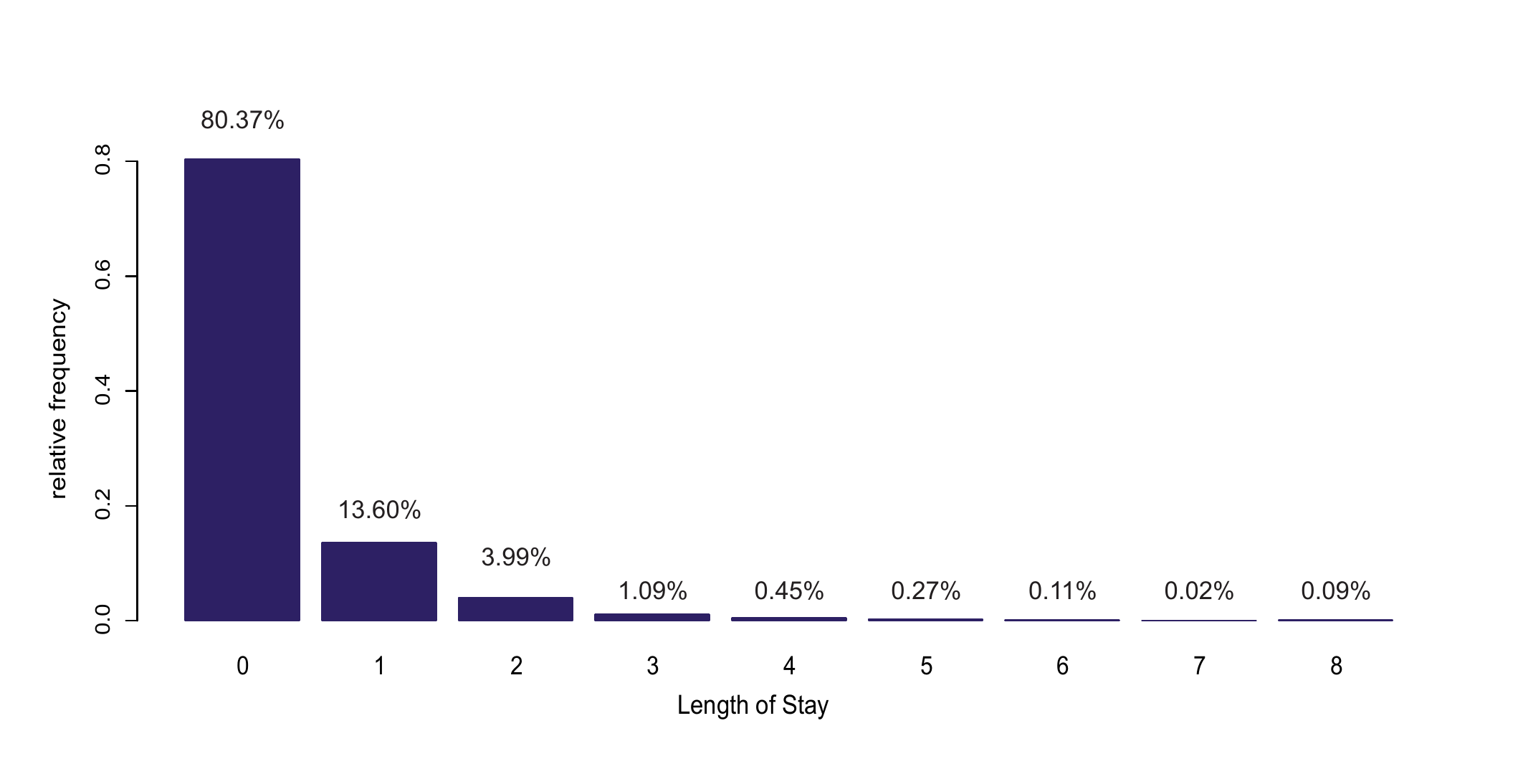}
		\captionof{figure}{Relative frequency of length of stay variable.}
	\end{center}
\end{figure}

\begin{table}[htp]
	\centering
	\caption{Gender-wise and health status-wise of HOSP variable} 
\begin{tabular}{lccccccc} \hline  \hline 
Gender	&$n$	& $\max$ &	$\min$	& mean	& variance	& ID	& \% of Zero\\ \hline 
Male	&1778	&8	     &0	        &0.311	& 0.554	    & 1.781	& 0.790 \\
Female	&2628	&8	    &0	        &0.286	& 0.5589	& 1.956	& 0.813 \\ \hline 
&&&&&&&\\
Health Status	&$n$	&$\max$	&$\min$	&mean &	variance& 	ID	& \% of Zero \\ \hline 
Poor	&554	&8	&0	&0.691	&1.389	&2.010	&0.610 \\
else	&3852	&8	&0	&0.239	&0.412	&1.723	&0.832 \\ \hline \hline 
\end{tabular}	
\end{table}

We fit  the proposed $\mathcal{UNB}$ model along with two competitive $\mathcal{UP}$ and $\mathcal{NB}$ regression models. The parameters are estimated by maximizing log-likelihood function (\ref{llc}) using \texttt{FindMaximum} function in Mathematica software package v.9.0 and the standard errors of the MLEs have been approximated by inverting the Hessian matrix evaluated at the MLEs. These result are given in Table 3. Further, we test $\mathscr{H}_0: \beta_k = 0$ against $\mathscr{H}_1: \beta_k \ne 0$, $k=0,1,\cdots,s$. We determine the significance of regression coefficients for all the three models by computing Wald's $t$-ratio, $\hat{\beta}_k/SE(\hat{\beta}_k)$, and its corresponding $p$-value at 5\% level of significance. Based on $p$-values, explanatory variables such as MARRIED, FAMINC, EMPLOYED, PRIVINS and MADICAID are not significant for $\mathcal{NB}$ and $\mathcal{UNB}$ regression models, whereas for $\mathcal{UP}$ models, variables MARRIED, FAMINC, EMPLOYED are not significant. The highest value of log-likelihood and lowest Akaike Information Criterion (AIC) value for $\mathcal{UNB}$ regression model confirm the good fit than other two models.

\begin{landscape}
 \begin{table}[]
 	\caption{Parameter estimates for Length of Stay in hospitalization data set under three regression models}
   \begin{tabular}{llrrrlrrrlrrr} \hline 
	\multicolumn{1}{c}{Models}     & \multicolumn{1}{c}{}            & \multicolumn{3}{c}{Uniform-Poisson $(\mathcal{UP})$}                                                           & \multicolumn{1}{c}{} & \multicolumn{3}{c}{Negative Binomial($\mathcal{NB}$)} & \multicolumn{1}{c}{} & \multicolumn{3}{c}{Uniform Negative Binomial ($\mathcal{UNB}$)}  \\ \hline 
 \multicolumn{1}{c}{Covariates} & \multicolumn{1}{c}{} & \multicolumn{1}{c}{Estimate(S.E.)} & \multicolumn{1}{c}{$t$-Wald} & \multicolumn{1}{c}{$p$-value} & \multicolumn{1}{c}{} & \multicolumn{1}{c}{Estimate(S.E.)} & \multicolumn{1}{c}{$t$-Wald} & \multicolumn{1}{c}{$p$-value} & \multicolumn{1}{c}{} & \multicolumn{1}{c}{Estimate(S.E.)} & \multicolumn{1}{c}{$t$-Wald} & \multicolumn{1}{c}{$p$-value} \\ \hline 
 & $\beta_0$                           & -3.530 (0.37)                      & -9.46                      & 0.000                           &                      & -3.795 (0.455)                     & -8.338                     & 0.000                       &                      & -3.811 (0.458)                     & 0.458                      & 0.000                           \\
EXCLHLTH                       & $\beta_1$                           & -0.725 (0.18)                      & -4.04                      & 0.000                           &                      & -0.703 (0.194)                     & -3.625                     & 0.000                       &                      & - 0.699 (0.195)                     & 0.195                      & 0.000                           \\
POORHLTH                       & $\beta_2$                           & 0.627 (0.07)                        & 8.44                       & 0.000                           &                      & 0.607 (0.095)                      & 6.401                      & 0.000                       &                      & 0.615 (0.096)                      & 0.096                      & 0.000                           \\
NUMCHRON                       & $\beta_3$                           & 0.274 (0.02)                       & 13.42                      & 0.000                           &                      & 0.286 (0.026)                      & 10.834                     & 0.000                       &                      & 0.289 (0.027)                      & 0.027                      & 0.000                           \\
AGE                            & $\beta_4$                           & 0.197 (0.04)                       & 4.2                        & 0.000                           &                      & 0.232 (0.57)                       & 4.034                      & 0.000                       &                      & 0.233 (0.058)                      & 0.058                      & 0.000                           \\
MALE                           & $\beta_5$                           & 0.154 (0.06)                       & 2.29                       & 0.02                        &                      & 0.187 (0.081)                      & 2.293                      & 0.022                       &                      & 0.188 (0.082)                      & 0.082                      & 0.022                       \\
MARRIED                        & $\beta_6$                           & -0.043 (0.07)                      & -0.62                      & 0.53                        &                      & -0.047 (0.085)                     & -0.559                     & 0.576                       &                      & -0.048 (0.086)                     & 0.086                      & 0.574                       \\
FAMINC                         & $\beta_7$                           & 0.005 (0.01)                       & 0.5                        & 0.61                        &                      & 0.002 (0.013)                      & 0.121                      & 0.903                       &                      & 0.002  (0.86)                      & 0.013                      & 0.891                       \\
EMPLOYED                       & $\beta_8$                           & 0.023 (0.11)                       & 0.2                        & 0.83                        &                      & 0.030 (0.130)                      & 0.229                      & 0.819                       &                      & 0.031 (0.132)                      & 0.132                      & 0.814                       \\
PRIVINS                        & $\beta_9$                           & 0.2 (0.08)                         & 2.38                       & 0.02                        &                      & 0.172 (0.102)                      & 1.685                      & 0.092                       &                      & 0.17 (0.102)                       & 0.102                      & 0.097                       \\
MEDICAID                       & $\beta_10$                          & 0.227 (0.11)                       & 2.08                       & 0.03                        &                      & 0.205 (0.136     )                      & 1.501                      & 0.133                       &                 & 0.208 (90.137)                      & 0.137                      & 0.130                        \\
	& $r$                               &                                    &                            &                             &                      & 1.755 (0.160)                      & 10.967                     & 0.000                       &                      & 0.884 (0.102)                      & 8.612                      & 0.000                           \\
&                                 &                                    &                            &                             &                      &                                    &                            &                             &                      &                                    &                            &                             \\
LL                             &                                 &                                    &                            & -2951.33                    &                      &                                    &                            & -2855.24                    &                      &                                    &                            & -2853.47                    \\
AIC                            &                                 &                                    &                            & 5924.660                    &                      &                                    &                            & 5734.48                     &                      &                                    &                            & 5730.94                     \\  \hline 
					\end{tabular}
				\end{table}
			\end{landscape}
\noindent For testing the closeness of $\mathcal{UNB}$ regression model with other competitive models, we use likelihood ratio test proposed by Voung (1989) with test statistics given as
	
\begin{equation} 
	\mathcal{Z}= \frac{1}{\omega \sqrt{n}}\left({l}_\mathcal{UNB}(\widehat{\Theta}_1)-{l}_{g}(\widehat{\Theta}_2)\right)  
\end{equation} 
	where 
	\[
	\omega^2=\frac{1}{n} \sum_{i=1}^{n} \left(\log\left(\frac{p(x_i|\widehat{\Theta}_1)}{g(x_i|\widehat{\Theta}_2)} \right)\right)^2 - \left(\frac{1}{n}\sum_{i=1}^{n}\log\left(\frac{p(x_i|\widehat{\Theta}_1)}{g(x_i|\widehat{\Theta}_2)} \right) \right)^2 
	\] 
	and $p$ and $g$ represent the PMF of $\mathcal{UNB}$ and other competitive model, respectively. Under null hypothesis, $\mathscr{H}_0: \mathbb{E}\left({l}_\mathcal{UNB}(\widehat{\Theta}_1)-{l}_{g}(\widehat{\Theta}_2)\right)= 0,$ against alternative hypothesis  $\mathscr{H}_1: \mathbb{E}\left({l}_\mathcal{UNB}(\widehat{\Theta}_1)-{l}_{g}(\widehat{\Theta}_2)\right)\ne 0$. Further $\mathcal{Z}$ is asymptotically normal distributed. The Voung test statistic value with corresponding $p$-values are given in Table 4. In both the cases $p$-value are less than 0.05 (significance level at 5\%). Hence we can strongly conclude that the proposed $\mathcal{UNB}$ regression model is preferred over other two competitive models.

\begin{table}[htp]
	\centering 
	\caption{Voung test results}  
	\begin{tabular}{ccc} \hline 
		& $\mathcal{Z}$ & $p$-value\\ \hline 
		$\mathcal{UNB}$ vs $\mathcal{NB}$ & 2.543 &  0.010 \\
		$\mathcal{UNB}$ vs $\mathcal{UP}$ & 4.887 & < 0.001  \\ \hline 
	\end{tabular}
\end{table}

\section{Final Comments} 
In this article we introduce a new discrete distribution for suitable for over dispersed and zero vertex count data. The new distribution is indexed by two parameters and possesses closed form expression for the probability generating function, moment generating function, mean and variance etc. Additionally, a new regression model in which the response variable follows the proposed count distribution is discussed and the results are compared with recently developed Uniform Poisson and negative binomial regression models.


\begin{thebibliography}{}

\bibitem{} Abramowitz, M. and  Stegun, I. (1972). \textit{Handbook of Mathematical Function}, 2nd ed., Dover, New York.

\bibitem{} 	Akdo$\breve{g}$an, Y., Ku$s$,C.,  Asgharzadeh, A., KInacI, I., \& Sharafi, F. (2016). Uniform-geometric distribution. \textit{Journal of Statistical Computation and Simulation}, 86(9), 1754-1770.

\bibitem{}	Ancarani, L.U. and Gasaneo, G. (2009). Derivatives of any order of the Gaussian Hypergeometric
function ${}_2F_1(a; b; c; z)$ with respect to the parameters $a$, $b$ and $c$, \textit{ Journal of Physics A: Mathematical and Theoretical}, 42, 395208-395218.

\bibitem{} Appell, P. and Kamp\'e de F\'eriet J. (1926). Fonctions Hyperg\'eom\'etriques et Hypersph\'eriques 	Polynomes d'Hermite (Paris: Gauthier-Villars).

\bibitem{} 
Berki, S.E., Ashcraft, M., and Newbrander, W.C. (1984). Length of siay variations within ICDA -8 diagnosis related groups, \textit{Medical Care}, 22, 126--142. 

\bibitem{} Bliss, C. I. and Fisher, R. A. (1953). Fitting the negative binomial distribution to biological data. \textit{Biometrics}, 9, 176--200.


\bibitem{}  Cameron, A. C. and Trivedi, P. K. (1998). \textit{Regression Analysis of Count Data}. Cambridge University Press, Cambridge.

\bibitem{} Chakraborty, S. and Ong, S. H. (2015) A COM-Poisson Generalization of the negative binomial Distribution.\textit{ Communications in Statistics-Theory and Methods}, 45, 4117-4135.

\bibitem{} Consul, P. C. (1989). \textit{Generalized Poisson Distribution: Properties and Applications}. Marcel Dekker, New York.

\bibitem{} Consul, P. C. and Famoye, F.  (1992). Generalized Poisson Regression Model, \textit{Communications in Statistics-Theory and Methods}, 21, 89-109.

\bibitem{} Deb, P. and Trivedi, P. K. (1997). Demand for medical care by the elderly: A finite mixture approach. \textit{Journal of  Applied Economics}, 12(3), 313–336.

\bibitem{} Ghitany, M. E., Al-Awadh, S. A. and Kalla, S. L. (2002) On Hypergeometric Generalized Negtive Binomial Distribution. \textit{International Journal of Mathematics and Mathematical Sciences}, 29, 727-736.

\bibitem{} G\'omez-D\'eniz, E. (2011). A new discrete distribution: properties and applications in medical care. \textit{ Journal of Applied Statistics}, 40(12), 2760-2770.

\bibitem{} Gupta, R. C. and Ong, S. H. (2004). A new generalization of the negative binomial distribution.\textit{Computational Statistics and Data Analysis}, 45, 287-300.


\bibitem{} Hilbe, J. M. (2011), Negative Binomial Regression, 2nd ed., Cambridge University Press.


\bibitem{} Jain, G.C. and Cosul P.C. (1971), A generalized negative binomial distribution. \textit{SIAM Journal of Applied Mathematics}, 21(4), 501-513.

\bibitem{} Karlis, D. and Xekalaki E. (2005), Mixed Poisson distributions.\textit{ International	Statistical Review}, 73, 35-58.



\bibitem{} Ong, S. H. and Lee P. A. (1986), On a generalized non-central negative binomial distribution. \textit{Communications in Statistics-Theory and Methods}, 15, 1065-1079.

\bibitem{} Vuong, Q. (1989). Likelihood ratio tests for model selection and non-nested hypotheses, \textit{Econometrica}, 57, 307-333.


\end{thebibliography}
\end{document}